\documentclass[11pt, reqno]{amsart}

 \usepackage{color}

\newcommand{\mysection}[1]{\section{#1}
      \setcounter{equation}{0}}

\newcommand{\nlimsup}{\operatornamewithlimits{\overline{lim}}}
\newcommand{\nliminf}{\operatornamewithlimits{\underline{lim}}}

\newtheorem{theorem}{Theorem}[section]
\newtheorem{lemma}[theorem]{Lemma}

\newtheorem{corollary}[theorem]{Corollary}

\theoremstyle{definition}
\newtheorem{assumption}{Assumption}[section]

\theoremstyle{remark}
\newtheorem{remark}{Remark}[section]

\newcommand\cbrk{\text{$]$\kern-.15em$]$}}
\newcommand\opar{\text{\raise.2ex\hbox{${\scriptstyle | }$}\kern-.34em$($} }

 \makeatletter
 \def\dashint{%
 \operatorname%
 {\,\,\text{\bf--}\kern-.98em\DOTSI\intop\ilimits@\!\!}}
 \makeatother

 \newcommand{\WO}{\overset{\scriptscriptstyle0}%
{ W}\,\!}

\newcommand\bR{\mathbb{R}}

\newcommand\cB{\mathcal{B}}

\newcommand\cF{\mathcal{F}}

\newcommand\cP{\mathcal{P}}

\newcommand\cN{\mathcal{N}}

\newcommand\dist{{\rm dist}\,}

\newcommand{\cWO}{\overset{\scriptscriptstyle0}%
{\mathcal{W}}\,\!}

\begin{document}

\title[Singularity of conditional distribution]{On singularity as a
function of time of a conditional distribution
of an exit time}

\author[N.V. Krylov]{N.V. Krylov}
\thanks{The  author was partially supported by
NSF Grant DNS-1160569}
\email{krylov@math.umn.edu}
\address{127 Vincent Hall, University of Minnesota,
 Minneapolis, MN, 55455}



\begin{abstract}
We establish Êthe singularity with
respect to Lebesgue measure
as a function of time of the conditional probability
that the sum of two one-dimensional Brownian motions will
exit from the unit interval before time $t$, given the trajectory
of the second Brownian motion up to the same time.
On the way of doing so we show that if one solves the one-dimensional
 heat
equation with zero condition on a trajectory of a one-dimensional
Brownian motion, which is the lateral boundary, then
for each moment of time with probability one
the normal derivative of the solution is zero,
provided that the diffusion of the 
Brownian motion is sufficiently
large.
\end{abstract}

\maketitle

\mysection{Main results}

Fix some constants $ \sigma,\sigma_{1}>0$   and consider the
equation
\begin{equation}
                                              \label{5.29.3}
 x_{t}=x_{0}+\sigma_{1}w_{t}+\sigma  b_{t},\quad
y_{t}=b_{t},
\end{equation}
  where $ w_{\cdot}$ and $b_{\cdot}$   are independent
one-dimensional standard Wiener processes,
$ x_{0}$ is  independent of the couple $ (w_{\cdot},b_{\cdot})$   and has
density $ \pi_{0}\in  C_{0}^{\infty}=C_{0}^{\infty}(G)$   
concentrated on $G$, where $ G=(0,1)$. Define
$$
 \cF^{b_{\cdot}}_{t}=\sigma(b_{s}:s\leq t),
\quad \tau=\inf\{t\geq0:x_{t}\not \in G\},
$$
  $$
 A_{t}=P(\tau\leq t\mid \cF^{b_{\cdot}}_{t}).
$$  

Here is our main result.
\begin{theorem}
                                             \label{theorem 5.10.01}
 There is a continuous and
nondecreasing modification of $ A_{t}$   which is singular
with respect to Lebesgue measure, the latter provided that
$ \sigma_{1}/\sigma$   is sufficiently small.
\end{theorem}

The process $A_{t}$ in a more general framework arose
in \cite{KW} as the main process governing
the conditional distribution of a signal process $x_{t}$
at the first time when it exists from a given domain.
In \cite{KW} the observations $y_{t}$ were only available
until the first exit time of $x_{t}$
from the domain. It turns out that in the setting of \eqref{5.29.3}
the conditional and the so-called 
unnormalized conditional distributions
of $x_{t}$ before it exits from $G$ given $y_{s},0\leq s \leq t$,
coincide. These unnormalized 
conditional distributions are known to satisfy
some linear stochastic partial differential equations
and then the properties of $A_{t}$ can be recovered
from some properties of solutions of these equations.

 To be more precise
for $(t,x)\in (0,\infty)\times(0,1)$ consider
the following (filtering) equation
\begin{equation}
                                              \label{3.17.1}
d\pi_{t}(x)=(1/2) a D^{2}\pi_{t}(x)\,dt
-\sigma D\pi_{t}(x)\,db_{t},
\end{equation}
where $a= \sigma_{1}^{2}+\sigma ^{2} $, with initial condition
$\pi_{0}(x)$ and zero lateral condition.

To explain in which sense we understand this equation, the initial
condition, and the boundary condition, we need some notation. 
Introduce the space 
$W^{1}_{2}=W^{1}_{2}(G)$ as the closure
of the set of continuously differentiable functions
in $\bar{G}$ in the norm
$$
\|u\|_{W^{1}_{2}}=\|u\|_{L_{2}}+\|Du\|_{L_{2}},
$$
where $Du$ is the derivative of $u$ and
$L_{2}=L_{2}(G)$, and we introduce
$\WO^{1}_{2}=\WO^{1}_{2}(G)$ as the
closure of $C^{\infty}_{0}=C^{\infty}_{0}(G)$ in the above norm.

Denote by $ \cP^{b_{\cdot}}$   the predictable
$ \sigma$-field in $ \Omega\times(0,\infty)$   associated
with the filtration $\{ \cF^{b_{\cdot}}_{t}\}$. For
$ T\in(0,\infty)$   introduce 
$$
 G_{T}=(0,T)\times G,\quad
\WO^{1}_{2}(G_{T})=L_{2}(\Omega\times(0,T),\cP^{b_{\cdot}},\WO^{1}_{2}).
$$
  We are looking for a function $ \pi_{t}(x)$   which is a
generalized function on $G$
 for each $ (\omega,t)$   such that $ \pi
\in\cap_{T}  \WO^{1}_{2}(G_{T})$     and
 for each $ \zeta\in C^{\infty}_{0} $   with probability one
for all $ t\in[0,\infty)$   it holds  that
\begin{equation}                                                  
                                            \label{1.16.1} 
(\pi_{t},\zeta)=(\pi_{0},\zeta)
-\int_{0}^{t}(1/2)( a  D \pi_{s} ,D
\zeta)\,ds-\int_{0}^{t} (\sigma D \pi_{s}
 ,\zeta)
\,db_{s},
\end{equation}
 where we use the notation
$$
(f,g)=\int_{G}f(x)g(x)\,dx.
$$

Observe that all expressions in \eqref{1.16.1} are well
defined due to the fact that the coefficients  of $ \pi$  
and of  
$ D \pi$   are constant and 
$$
 \pi, D \pi\in\L_{2}(T):= L_{2}(\Omega\times(0,T),
\cP^{b_{\cdot}},L_{2})
$$
  for any $ T\in(0,\infty)$.

Recall that by assumption $\pi_{0}\in C_{0}^{\infty}$.

\begin{theorem}
                                           \label{theorem 3.19.1}

In the class $\bigcap_{T}\cWO^{1}_{2}(G_{T})$ there exists a unique solution
$\pi_{t}$ of equation \eqref{3.17.1} with initial condition
$\pi_{0}$. In addition,   
$
\pi_{t}\geq0
$
for all $t\in[0,\infty)$ (a.s.). With probability
one $\pi_{t}$ is continuous in $L_{1}=L_{1}(G)$
and in $L_{2}$.
\end{theorem}
The existence, uniqueness, and the (a.s.)
continuity in
$L_{2}$ of $\pi$ is a classical result proved in
many places in a variety of settings (see, for
instance,
\cite{Pa1}, \cite{KR}, \cite{Ro}, and the
references therein). That $\pi_{t}$ is (a.s.)
continuous as an $L_{1}$-function follows from its
$L_{2}$-continuity and the boundedness of $G$. The
fact that $\pi\geq0$ follows from the maximum
principle (see, for instance, Theorem 1.1 of
\cite{Kr07}) and the fact that, if $u\in
\WO^{1}_{2}$, then $u^{+}\in \WO^{1}_{2}$.  

The connection of $A_{t}$ to $\pi_{t}$
is established on the basis of Lemma 1.8 of \cite{KW},
which   our situation
reads as follows.

\begin{lemma}
                                                 \label{lemma 10.11.2}
For any Borel bounded or nonnegative function 
$\phi$ on $G$ and $t\in[0,\infty)$ we have (a.s.)
\begin{equation}
                                                       \label{10.11.4}
E\big\{I_{\tau>t}\phi(x_{t})
\mid \cF^{b_{\cdot}}_{t}\}=(\pi_{t},\phi).
\end{equation}
In particular,   for each $t\in[0,\infty)$ (a.s.)
\begin{equation}
                                                       \label{10.11.6}
P\{\tau>t\mid\cF^{b_{\cdot}}_{t}\}
=(\pi_{t},1) .
\end{equation}
Finally, (a.s.) 
we have
$ (\pi_{t},1)>0$ for all $t\in[0,\infty)$.
 
\end{lemma}

By Lemma \ref{lemma 10.11.2} for any $t\in[0,\infty)$
$$
P(\tau\leq t\mid \cF^{b_{\cdot}}_{t})=1-(\pi_{t},1)
$$
(a.s.) and by Theorem \ref{theorem 3.19.1} the right-hand side
is continuous in $t$ (a.s). Also
it turns out (see \cite{KW}) that  the process
$ (\pi_{t},1)$ is decreasing (a.s). Therefore,
in Theorem \ref{theorem 5.10.01} by the modification
of $A_{t}$, which we identify with the original $A_{t}$,
we mean $1-(\pi_{t},1)$.

Observe that if in \eqref{1.16.1}
we were allowed to first integrate by parts
to replace 
$$
(D\pi_{s},D\zeta)\quad\text{with}\quad-(D^{2}\pi_{s},\zeta)
$$
and then in the so modified version of \eqref{1.16.1}
take $\zeta\equiv1$
($\not\in C^{\infty}_{0}$), then we would formally
obtain that
\begin{equation}
                                                 \label{5.29.4}
A_{t}=1-(\pi_{t},1)=(1/2)a\int_{0}^{t}[D\pi_{s}(0)-
D\pi_{s}(1)]\,ds.
\end{equation}

This shows that $A_{t}$ is related to the normal
derivative of $\pi_{s}$ on the boundary of $G$,
investigating which is done on the basis of a
  different description of $\pi_{s}$.

We are going to state our second main result, which
is about solutions of the heat equation in
curvilinear cylinders whose lateral boundary is a trajectory
of a Wiener process.

\begin{theorem}
                                  \label{theorem 3.21.1}
For almost any $\omega$ there exists
a unique function $u_{t}(x)$ defined, bounded,  and continuous
in the closure of
$$
\Gamma(b_{\cdot})=
\{(t,x):t\geq 0,\quad x\in 
(\sigma b_{t},1+\sigma b_{t})\}
$$ 
such that it is infinitely differentiable
with respect to $(t,x)$ in $\Gamma(b_{\cdot})$, satisfies
there the equation
\begin{equation}
                                      \label{3.21.2}
\partial_{t}u_{t}=(1/2)\sigma_{1}^{2}D^{2}u
\end{equation}
and satisfies the conditions  $u_{0}(x)=\pi_{0}
(x)$, $x\in[0,1]$, and $u_{t}(\sigma b_{t})
=u_{t}(1+\sigma b_{t})=0$, $t\geq0$.
Furthermore, if $\sigma_{1}/\sigma$ is sufficiently
small, then for any $t\in[0,\infty)$
\begin{equation}
                                       \label{3.22.1}
\lim_{x\downarrow \sigma b_{t}}\frac{u_{t}(x)}{x-\sigma b_{t}}
=\lim_{x\uparrow 1+\sigma b_{t}}\frac{u_{t}(x)}{
1+\sigma b_{t}-x }=0
\end{equation}
almost surely, so that the derivative of $u_{t}(x)$
on the boundary of $\Gamma(b_{\cdot})$ is zero
(a.s.) for any fixed $t$.

Finally,
with probability one
\begin{equation}
                                        \label{6.14.1}
\int_{G}|u_{t}(x-\sigma b_{t})-\pi_{t}(x)|^{2}\,dx=0
\end{equation}
for all $t\geq0$, so that $u_{t}(x-\sigma b_{t})$
is a modification of $\pi_{t}(x)$, and for this modification,
for any $t\in[0,\infty)$,
\begin{equation}
                                       \label{5.29.6}
\lim_{x\downarrow0}\frac{\pi_{t}(x)}{x}
=\lim_{x\uparrow1}\frac{\pi_{t}(x)}{1-x}=0
\end{equation}
almost surely if $\sigma_{1}/\sigma$ is sufficiently
small.

\end{theorem}

The last statement of Theorem \ref{theorem 3.21.1}
makes the representation \eqref{5.29.4} questionable,
and, even though there is a limit procedure  showing that
\eqref{5.29.4} holds in a generalized sense similar to that of
the local time for Brownian motion (see \cite{KW}),
one would rightfully suspect that $A_{t}$ is not absolutely
continuous with respect to $t$.

The last statement of Theorem \ref{theorem 3.21.1}
should not make the reader over optimistic about the 
continuity properties of $\pi_{t}(x)$ in $x$ near the boundary
of $G$ (see Remark \ref{remark 3.22.1}).

Still the following theorem will be easily derived from
known results. 
Take $\alpha\in(0,1)$  and $c\in(0,\infty)$ and introduce
$$
p=p(c):=P(\sup_{t\leq1}w_{t}-
\inf_{t\leq1}w_{t}\geq c(\sqrt{2}-1)/2),
\quad r=r(\alpha,c):=\frac{\alpha(1-p )}{p (1-\alpha)},
$$
$$
\beta=\beta(\alpha,c):=2\frac{(r-1)p +1}{r^{\alpha}}
=2\frac{1-p}{1-\alpha}r^{-\alpha}.
$$
As is  easy to see, for any $\alpha\in(0,1)$,
we have $p(c)\to0$, $r(\alpha,c)\to\infty$,
and $\beta(\alpha,c)\to0$
as $c\to\infty$. 
It follows from \cite{Kr_04} that there exists a function
  $\alpha(c)$,
$c\in(0,\infty)$, with values in $(0,1)$
such that $\alpha(c)\to0$ as $c\to\infty$ and $\alpha(c)\leq
\alpha$ for any $\alpha$ satisfying 
$\beta(\alpha ,c)<1$.

Next, take some constants $c\geq0$, $d>0$ and
for $x\in\bR$ define
$$
\tau_{d,x}=\inf\{t\geq0:d/\sqrt{2}+w_{t}=x \},
$$
$$
\gamma(c,d,\delta)=P(\tau_{d,d} \wedge(\delta/2)<\tau_{d,-c} ).
$$
 
\begin{theorem}
                                  \label{theorem 3.21.01}
Take the modification of $\pi_{t}$
from Theorem \ref{theorem 3.21.1}. Then 
 for any $t\in[0,\infty)$, $c,d>0$, such that 
$\alpha(c \varepsilon )<1$, where $\varepsilon=\sigma_{1}/\sigma$,
and
$\nu$ satisfying
$$
0\leq\nu< \nu_{0}:=
(1-\alpha(c \varepsilon))\log_{2}\gamma^{-2}(c,d,1),
$$
we have that with probability one
\begin{equation}
                                         \label{3.24.1}
\sup_{x\in(0,1)}\sup_{t\in[0,T]}\frac{\pi_{t}( x)}
{x^{\nu}}<\infty.
\end{equation}

Furthermore,  there exists a
 constant  $ \nu\in(0,\infty)$ such that for any $T\in(0,\infty)$  

\begin{equation}
                                           \label{3.24.02} 
 E\sup_{x\in(0,1)}\sup_{t\in[0,T]}\frac{\pi_{t}(x)}
{x^{\nu}(1-x)^{\nu}}<\infty.
\end{equation}

\end{theorem}

\begin{remark}
                             \label{remark 3.22.1}
The largest possible value  of $\nu$ in
(\ref{3.24.1}) is  unknown. However,
Theorem 5.1 and Lemma 4.1 of \cite{Kr03}
show that if we take a $\mu>0$
and
$$
\nu=(1+\mu)(2\pi\varepsilon^{2})^{-1/2}
e^{-1/(2\varepsilon^{2})},
$$
then for $\varepsilon=\sigma_{1}/\sigma$ small enough the left-hand side
of (\ref{3.24.1}) equals infinity with probability one.
Therefore, the largest value of $\nu$
is extremely small if $\varepsilon$ is small.
\end{remark}

\begin{remark}
                             \label{remark 3.22.2}
The fact that equation \eqref{3.22.1}
holds (a.s.) does not contradict Remark
\ref{remark 3.22.1}, because \eqref{3.24.1}
gives
 an estimate which is uniform with respect to $t$,
and on almost each trajectory of $b_{\cdot}$ there are
points $t$ such that $v_{t}(x)/x\to\infty$
as $x\downarrow 0$.
\end{remark}

We prove Theorem \ref{theorem 5.10.01}
in Section \ref{section 6.14.1}
assuming that Theorems \ref{theorem 3.21.1}
and \ref{theorem 3.21.01}
are true. In Section \ref{section 6.14.2}
we prove Theorems \ref{theorem 3.21.1} and
\ref{theorem 3.21.01}. 
The first assertion of Theorem  \ref{theorem 3.21.1}
and It\^o's formula easily lead to the conclusion that 
$u_{t}(x-\sigma b_{t})$ is a classical solution
of \eqref{3.17.1} and the
assertion concerning \eqref{6.14.1}
is proved by showing that classical solutions
coincide with generalized ones in a much more general situation
in Section \ref{section 6.14.3}.

\mysection{Proof of Theorem \protect\ref{theorem 5.10.01}}
                                    \label{section 6.14.1}

 We start by proving that for each $t_{0}\in(0,\infty)$
and $t_{n}=t_{0}+1/n$
with probability one
\begin{equation}
                                          \label{3.22.3}
E\big(\nliminf_{n\to\infty}\frac{1}{t_{n}-t_{0}}(A_{t_{n}}-A_{t_{0}}) 
\mid \cF^{b_{\cdot}}_{t_{0}}\big)=0.
\end{equation}

Observe that for any 
$\zeta\in C^{\infty}_{0}$ and $t>t_{0}$ 
\begin{equation}
                                    \label{3.22.4}
(\pi_{t},\zeta)=(\pi_{t_{0}},\zeta)
+\int_{t_{0}}^{t}(
 \pi_{s}, (1/2)aD^{2}\zeta)\,ds+\int_{t_{0}}^{t}
( \pi_{s},\sigma  D \zeta)
\,db_{s},
\end{equation}
where and below we are dealing with the modification of $\pi_{t}$
from Theorem \ref{theorem 3.21.1}.
We multiply both parts of this equation by the 
indicator function of a set $F\in\cF^{b_{\cdot}}_{t_{0}}$
and then take the   expectations of
both parts. Then by denoting
$$
\phi^{F}_{t}(x)=E\pi_{t}(x)I_{F}
$$
we find
\begin{equation}
                                    \label{3.22.5}
(\phi^{F}_{t},\zeta)=(\phi^{F}_{t_{0}},\zeta)
+\int_{t_{0}}^{t}(
\phi^{F}_{s}, (1/2)aD^{2}\zeta)\,ds.
\end{equation}
Observe that $\phi^{F}_{t}(x)$ is continuous
in $\bar{G}_{\infty}$ because such is $\pi$
which is in addition bounded.
Estimate \eqref{3.24.02}  shows that 
$\phi^{F}_{t}(x)\to0$ as $x\to\{0,1\}$, $x\in(0,1)$,
$t\geq0$. Thus $\phi^{F}_{t}$ is
 a continuous in $[t_{0},\infty)\times\bar{G}$
weak solution of the equation
\begin{equation}
                                    \label{3.22.6}
\partial_{t}\eta_{t}=(1/2)aD^{2}\eta_{t}.
\end{equation}
By uniqueness of  such solutions, $\phi^{F}_{t}$
is a classical solution of this equation with
zero boundary data.

By the maximum principle we have
$\phi^{F}_{t}\leq \psi^{F}_{t}$, $t\geq t_{0}$,
where $\psi^{F}_{t}$ is defined as a unique bounded classical
solution of \eqref{3.22.6} for $t\geq t_{0}$, $x>0$,
with initial data $\psi^{F}_{t_{0}}(x)=\phi^{F}_{t_{0}}(x)
I_{(0,1)}(x)$ and zero boundary condition.

The following explicit representation for such solutions
is well known:
$$
\psi^{F}_{t}(x)=\frac{1}{\sqrt{2\pi a(t-t_{0})}}
\int_{0}^{1}
\phi^{F}_{t_{0}}(y) \exp\big[-\frac{(x-y)^{2}}{2 a(t-t_{0})}\big]
-\exp\big[-\frac{(x+y)^{2}}{2 a (t-t_{0})}\big]\,dy.
$$

Furthermore,
$$
E(A_{t}-A_{t_{0}})I_{F}=
\int_{0}^{1}[\phi^{F}_{t}(x)-
\phi^{F}_{t_{0}}(x)]\,dx
\leq \int_{0}^{\infty}[\psi^{F}_{t}(x)-
\phi^{F}_{t_{0}}(x)I_{(0,1)}(x)]\,dx.
$$
Observe that
$$
\partial_{t}\int_{0}^{\infty}[\psi^{F}_{t}(x)-
\phi^{F}_{t_{0}}(x)I_{(0,1)}(x)]\,dx=
(1/2)a\int_{0}^{\infty}D^{2}\psi^{F}_{t}(x)\,dx
$$
$$
=(1/2)aD\psi^{F}_{t}(0),
$$
so that
$$
E(A_{t}-A_{t_{0}})I_{F}\leq  \int_{t_{0}}^{t}
\frac{1}{\sqrt{2\pi a(r-t_{0})^{3}}}  
\int_{0}^{1}y\phi^{F}_{t_{0}}(y)e^{-y^{2}/(2ar-2at_{0})}\,dy dr
$$
$$
=\int_{t_{0}}^{t}
\frac{2a}{\sqrt{2\pi a(r-t_{0})}}  
\int_{0}^{1/\sqrt{2ar-2at_{0}}}
\phi^{F}_{t_{0}}(x\sqrt{2ar-2at_{0}})xe^{-x^{2}}\,dxdr,
$$
which after taking into account the arbitrariness
of $F\in\cF^{b_{\cdot}}_{t_{0}}$ leads to
$$
E(A_{t}-A_{t_{0}}\mid \cF^{b_{\cdot}}_{t_{0}})
$$
$$ \leq
\int_{t_{0}}^{t}
\frac{2a}{\sqrt{2\pi a(r-t_{0})}}  
\int_{0}^{\infty}
\pi_{t_{0}}(x\sqrt{2ar-2at_{0}})I_{x\sqrt{2ar-2at_{0}}\leq1}xe^{-x^{2}}\,dxdr
$$
almost surely for any $t>t_{0}$. By Theorem \ref{theorem 3.21.1} with probability one
$\pi_{t_{0}}(x)=x\theta (x)$, $x\in[0,1]$, where $\theta $ is a bounded function
of $x$ tending to zero as $x\downarrow0$. Hence
$$
E(A_{t}-A_{t_{0}}\mid \cF^{b_{\cdot}}_{t_{0}})
$$
$$ \leq \frac{2a}{\sqrt{ \pi }} 
\int_{t_{0}}^{t}
 \int_{0}^{\infty}\theta(\sqrt{2arx-2at_{0}x}\wedge 1)x^{2}e^{-x^{2}}
\,dxdr .
$$

By the dominated convergence theorem (a.s.)
$$
\lim_{r\downarrow t_{0}}
 \int_{0}^{\infty}\theta(\sqrt{2arx-2at_{0}x}\wedge 1)x^{2}e^{-x^{2}}
\,dx=0
$$
implying that (a.s.)
$$
\lim_{n\to\infty}\frac{1}{t_{n}-t_{0}}\int_{t_{0}}^{t_{n}}
 \int_{0}^{\infty}\theta(\sqrt{2arx-2at_{0}x}\wedge 1)x^{2}e^{-x^{2}}
\,dxdr,
$$
$$
\nlimsup_{n\to\infty}\frac{1}{t_{n}-t_{0}}
E(A_{t_{n}}-A_{t_{0}}\mid \cF^{b_{\cdot}}_{t_{0}})=0,
$$
which yields \eqref{3.22.3} by Fatou's lemma.

Thus, for any $t\geq0$, for almost all $\omega$
\begin{equation}
                                                \label{3.16.4}
\nliminf_{n\to\infty}n(A_{t+1/n}-A_{t})=0.
\end{equation}
By Fubini's theorem, for almost any $\omega$, equation \eqref{3.16.4}
holds for almost all $t$. It follows that the derivative of
$A_{t}$ is zero for almost all $t$ and the theorem is proved.

\mysection{Proof of Theorems  
\protect\ref{theorem 3.21.1} and \protect\ref{theorem 3.21.01}}
                                    \label{section 6.14.2}
 
{\bf Proof of Theorem \ref{theorem 3.21.1}}.
On the space $C$ of continuous functions on
$[0,\infty)$ with   Wiener measure $W$
introduce the coordinate process $x_{t}(x_{\cdot}):=x_{t}$,
which is a Wiener process. 
For $t\geq0$, $x\in\bR$, and $x_{\cdot},y_{\cdot}\in C$
such that $y_{0}=0$ define
$$
\tau(t,x,x_{\cdot},y_{\cdot})=\inf\{s\geq0:x+
\sigma_{1}x_{s}\not\in( \sigma y_{t-s},
\sigma y_{t-s}+1)\},
$$
where $y_{r}:=y_{0}$ for $r\leq0$.
Then the function
$$
u_{t}(y_{\cdot} , x):=
\int_{C}\pi_{0}(x_{t})I_{ \tau(t,x,x_{\cdot},y_{\cdot} )\geq t}
\,W(dx_{\cdot}).
$$
is the probabilistic solution of the heat equation
$$
\partial_{t}u_{t}=(1/2)\sigma_{1}^{2}D^{2}u_{t}
$$
in $\Gamma(y_{\cdot})$ 
with boundary conditions
$$
u_{t}( y_{t})=  u_{t}( y_{t}+1)=0,\quad t>0,
$$
$$
u_{0}( x)=\pi_{0}(x),\quad 0\leq x\leq 1.
$$
 Due to interior
estimates of derivatives of solutions to the heat equation,
$u_{t}(y_{\cdot} , x)$ is infinitely differentiable in $\Gamma(y_{\cdot})$.
Its continuity up to $\{0\}\times[0,1]$ easily follows from
the fact that $\pi_{0}$ is a continuous function on $[0,1]$ vanishing at 
$0$ and $1$. The continuity
of its derivatives up to $\{0\}\times(0,1)$
follows from the fact that $\pi_{0}$ is infinitely
differentiable. Next, as in the proof of Theorem 4.1 of \cite{Kr_04}
one shows that for constants $\nu$ 
(different in \eqref{6.12.1} and \eqref{6.12.2})
as in Theorem \ref{theorem 3.21.01}
we have that,
 for almost any trajectory of $b_{\cdot}$,
\begin{equation}
                                                     \label{6.12.1}
\sup_{x\in(0,1)}\sup_{t\in[0,T]}\frac{u_{t}(b_{\cdot} , x)}
{(b_{t}+x)^{\nu}(b_{t}+1-x)^{\nu}}<\infty ,
\end{equation}
\begin{equation}
                                                     \label{6.12.2}
E\sup_{x\in(0,1)}\sup_{t\in[0,T]}\frac{u_{t}(b_{\cdot} , x)}
{(b_{t}+x)^{\nu}(b_{t}+1-x)^{\nu}}<\infty .
\end{equation}
In particular, with probability one
$u_{t}(b_{\cdot} , x)$ is continuous at the lateral boundary
of $\Gamma(b_{\cdot})$.  

Next we deal with \eqref{3.22.1}
We organize the proof in the following way.
For $t\geq0$, $x\in\bR$, and $x_{\cdot},y_{\cdot}\in C$
such that $y_{0}=0$ define
$$
\gamma(t,x,x_{\cdot},y_{\cdot})=\inf\{s\geq0:x+
\sigma_{1}x_{s}\leq\sigma y_{t-s}\},
$$
where $y_{r}:=y_{0}$ for $r\leq0$.
Also let
$$
v_{t}(y_{\cdot} , x):=
\int_{C}\pi_{0}I_{G}(x_{t})I_{ \tau(t,x,x_{\cdot},y_{\cdot} )\geq t}
\,W(dx_{\cdot}).
$$

\begin{lemma}
                                   \label{lemma 4.5.1}
Let $B_{t}$ be a one-dimensional Wiener process
and $\gamma\in(0,1)$.
Then with probability one there exists
a sequence of integers $0\leq m_{1}<m_{2}<...$ such that
$B_{t_{k}}\geq\sqrt{t_{k}}$ for all $k$, where $t_{k}=\gamma^{m_{k}}$.
Moreover, $m_{k}\leq\beta k$
for all sufficiently large $k$, where $\beta$ is any fixed number such 
that $\alpha\beta>1$,
$\alpha=P(B_{1}\geq1)$.
 \end{lemma}

Proof. The sequence $I_{B_{\gamma^{m}}\geq \gamma^{m/2}},m=0,1,...$,
is stationary, so that the limit
$$
\lim_{n\to\infty}\frac{1}{m}\sum_{k=1}^{m}I_{B_{\gamma^{k}}
\geq \gamma^{k/2}}
$$
exists (a.s.). By the 0-1 law this limit is a constant 
and equals $\alpha$ (a.s.). Set
$$
m_{1}=\inf\{k\geq1:B_{\gamma^{k}}\geq \gamma^{k/2}\},\quad
m_{n+1}=\inf\{k>m_{n}:B_{\gamma^{k}}\geq \gamma^{k/2}\}.
$$
Then the number of $k\in\{1,2,...\}$ such that $m_{k}\leq m$
divided by $m$ tends to $\alpha$. It follows that
the number of integers $i\in \{1,2,... ,\beta k\}$
such that $m_{i}\leq\beta k$
for all large $k$   is greater
than $\beta'k\alpha$, where $\beta'$ is any number
such that $\beta'<\beta$ and $\beta'\alpha>1$.
On the other hand
there are always exactly
 $k$ values of $i\in\{1,2,..., m_{k}\}$
such that $B_{\gamma^{i}}\geq \gamma^{i/2}$.
 Since $\beta'\alpha>1$, it follows that for any
sufficiently large $k$
the inequality $m_{k}\geq \beta k$ is impossible.
The lemma is proved.

\begin{lemma}
                                   \label{lemma 3.29.1}
If $\sigma_{1}/\sigma$ is sufficiently small, then for each
$T\in(0,\infty)$
\begin{equation}
                                       \label{3.29.3}
\lim_{x\downarrow0}\frac{v_{T}(x,b_{\cdot})}{x-\sigma b_{T}}=
0
\end{equation}
almost surely.
\end{lemma}

Proof. Fix a $T\in(0,\infty)$, set $\gamma=1/2$, and 
take integers $0\leq m_{1}<m_{2}<...$ such that
$b_{T-t_{k}}-b_{T}\geq\sqrt{t_{k}}$ for all $k$, 
where $t_{k}=\gamma^{m_{k}}$
and $m_{k}\leq\beta k$ for all large $k$.
By Lemma \ref{lemma 4.5.1}
such a sequence exists with probability one.
Then notice that the inequality 
$\tau(T,x,x_{\cdot},b_{\cdot})\geq T$ implies that
$$
x+\sigma_{1} x_{t_{k}}\geq \sigma b_{T-t_{k}}
\geq \sigma b_{T}+\sigma \sqrt{ t_{k}}
$$
 for all $k$ such that $t_{k}\leq T$. Denote by $k_{0}$
the smallest $k$ such that $t_{k}\leq T$.
We also take into account that $\pi_{0}$ is a bounded function
and conclude that for any integer $n\geq k_{0}$
$$
v_{T}(x,b_{\cdot})\leq NW(z+\sigma_{1} x_{t_{k}}-\sigma 
\sqrt{t_{k}}\geq0,k=k_{0},...,n),
$$
where $N=\sup\pi_{0}$ and $z=x-\sigma b_{T}$.
It follows that to prove \eqref{3.29.3} it suffices to show
that there exists an integer-valued
 function $n=n(x)\geq k_{0}$ such that
\begin{equation}
                                                     \label{3.29.5}
\lim_{x\downarrow0}\frac{1}{x}
P(x+  w_{t_{k}}- K
\sqrt{t_{k}}\geq0,k=k_{0},...,n(x))=0,
\end{equation}
where $K=\sigma/\sigma_{1}=\varepsilon^{-1}$.

 Observe that by Girsanov's theorem for any $n\geq k_{0}+1$
$$
P(x+  w_{t_{k}}- K
\sqrt{t_{k}}\geq0,k=k_{0},k_{0}+1,...,n )
$$
$$
= EI_{\Gamma_{n}(x)}
\exp\big(-\int_{t_{n}}^{t_{k_{0}}}f'(t)\,dw_{t}
-(1/2)\int_{t_{n}}^{t_{k_{0}}}[f'(t)]^{2}\,dt\big)=:I_{n}(x),
$$
where $f(t)= K\sqrt{t}$ for $t\geq t_{n}$, $f(t)= K\sqrt{t_{n}}$
for $t\in[0,t_{n}]$,   and 
$$
\Gamma_{n}(x)=\{\omega:x+  w_{t_{k}} 
- K\sqrt{t_{n}}  \geq 0,k=k_{0},...,n\}
$$
 
Next, note that for bounded nonrandom functions $g$
$$
E\{\exp\int_{s}^{t}g(r)\,dw_{r}
\mid w_{s},w_{t}\}
$$
$$
=\exp\bigg(\frac{w_{t}-w_{s}}{t-s}\int_{s}^{t}g(u)\,du
+(1/2)
\int_{s}^{t}\big[g(r)-\frac{1}{t-s}\int_{s}^{t}g(u)\,du\big]^{2}\,dr
\bigg )
$$
because
$$
\int_{s}^{t}\big[g(r)-\frac{1}{t-s}\int_{s}^{t}g(u)\,du\big]\,dw_{r}
$$
is independent of $w_{s},w_{t}$. Then we use the fact that
$$
\int_{s}^{t}\big[g(r)-\frac{1}{t-s}\int_{s}^{t}g(u)\,du\big]^{2}\,dr
=\int_{s}^{t}g^{2}(r)\,dr-\frac{1}{t-s}
\big(\int_{s}^{t}g(u)\,du\big)^{2}
$$
Hence, by applying this to   $g=-f'$ we find that
$$
I_{n}(x)=EI_{\Gamma_{n}(x)}
\exp\big(D_{n}-(1/2)C_{n})
$$
with
$$
D_{n}:=- K\sum_{k=k_{0}+1}^{n}s_{k}(w_{t_{k-1}}-w_{t_{k}}),\quad s_{k}=
(\sqrt{t_{k-1}}+\sqrt{t_{k}})^{-1},
$$
$$
C_{n}:=\sum_{k=k_{0}+1}^{n} \frac{1}{t_{k-1}-t_{k}}
\big(\int_{t_{k}}^{t_{k-1}}f'(u)\,du\big)^{2}
$$
$$
= K^{2}\sum_{k=k_{0}+1}^{n} \frac{\sqrt{t_{k-1}}-\sqrt{t_{k}}}
{\sqrt{t_{k-1}}+\sqrt{t_{k}}}
\geq  K^{2}\frac{1-\sqrt{\gamma}}{1+\sqrt{\gamma}}(n-k_{0})=
:2K\kappa (n-k_{0}),
$$
where the inequality follows from the fact that $\gamma t_{k-1}
\geq  t_{k}$ and $\kappa$ is defined by the last equality.

Now we use summation by parts to see that
$$
 D_{n}=- Kw_{t_{k_{0}}}s_{k_{0}+1}
+ Kw_{t_{n}}s_{n}
- K\sum_{k=k_{0}+1}^{n-1}w_{t_{k}} (s_{k+1}-s_{k}).
$$
On the event   $\Gamma_{n}(x)$ this quantity
is smaller than
$$
 Kw_{t_{n}}s_{n}+
  K(x-K\sqrt{t_{n}}) \big(s_{k_{0}+1}
+
\sum_{k=k_{0}+1}^{n-1}  (s_{k+1}-s_{k}) \big) 
$$
$$
= Kw_{t_{n}}s_{n}+ K(x-K\sqrt{t_{n}}) s_{n}.
$$

It follows that
$$
I_{n}(x)\leq EI_{\Gamma_{n}(x)}
\exp\big(   Kw_{t_{n}}s_{n}
+ Kx s_{n}
- K^{2}\kappa(n-k_{0})\big) ,
$$
$$
\leq E 
\exp\big(   Kw_{t_{n}}s_{n}
+ Kx s_{n}
- K^{2}\kappa(n-k_{0})\big)
$$
$$
= \exp\big( ( K^{2}/2)t_{n}s^{2}_{n}
+
 Kx s_{n}- K^{2}\kappa(n-k_{0})\big).
$$

Now it is time to choose $n=n(x)$. We take $n=n(x)$
so that $x^{2}\in[t_{n},t_{n-1}]$. Then $t_{n}s_{n}\leq\sqrt{t_{n}}
\leq x$ and $xs_{n}\leq x/\sqrt{t_{n-1}}\leq 1$. Also
$t_{n}s^{2}_{n}\leq1$ and 
since $t_{n}\geq \gamma^{\beta n}$, we have $x\geq \gamma^{\beta n/2}$,
which implies that
$$
\lim_{x\downarrow0}x^{-1}I_{n(x)}(x)
$$
$$
\leq
\lim_{n\to\infty}\exp\big( K^{2}/2+  K- K^{2}\kappa(n-k_{0})-
(1/2)\beta n\ln \gamma\big)=0
$$
if $\beta|\ln \gamma|< 2K^{2}\kappa$, which is true if $\sigma_{1}/\sigma$
is small enough.
This proves the lemma.

\begin{corollary}
                                        \label{corollary 6.13.1}

If $\sigma_{1}/\sigma$ is sufficiently small, then
\eqref{3.22.1} holds (a.s.) for any fixed $t\geq0$.
\end{corollary}

Indeed, the equality of the extreme terms in \eqref{3.22.1}
follows from Lemma \ref{lemma 3.29.1} since $u_{t}\leq v_{t}$.
The remaining equality is proved similarly by replacing $x$
with $1-x$.

It only remains to prove the last assertion
of the theorem.

Observe that
 $\tau(t,x,x_{\cdot},y_{\cdot})$
is a lower semicontinuous function
of its arguments. Therefore by Fubini's theorem
$u_{t}(y_{\cdot},x)$ is a Borel function of $(y_{\cdot},t,x)$.
Furthermore, $u_{t}(y_{\cdot}, x)$ will not change if
we change $y_{r}$ for $r>t$. Hence,
$u_{t}(y_{\cdot}, x)$ is $\cN_{t}$-measurable,
where $\cN_{t}=\sigma(y_{r}:r\leq t,y_{\cdot}\in C)$.
Therefore, 
$$
v_{t}(x):=u_{t}(b_{\cdot},x-\sigma b_{t})
$$
 is $\cF_{t}$-measurable for each 
$(t,x)\in\bar{G}_{\infty}$. After that 
 in the same way the usual It\^o's formula is proved
on the basis of Taylor's formula and the fact that $u_{t}(y_{\cdot},x)$
is infinitely differentiable we obtain that for any
$x\in G$ almost surely for all $t\geq0$
$$
v_{t}(x)=\pi_{0}(x)+(a/2)\int_{0}^{t}D^{2}v_{s}(x)\,ds
-\sigma\int_{0}^{t}D v_{s}(x)\,db_{s}.
$$
 The above properties of $u_{t}$ and Theorem \ref{theorem 5.28.1}
now imply that (perhaps after modifying $v_{t}$ on a set of probability 
zero) $v_{t}$ satisfies \eqref{3.17.1}
with zero boundary condition 
and initial condition $\pi_{0}$ in the sense explained below
that formula and is such that
$$
\int_{0}^{T}\|v_{t}\|^{2}_{W^{1}_{2}(G)}\,dt<\infty
$$
 for each $T\geq0$. Uniqueness of solutions
of \eqref{3.17.1} in this class of functions is
a classical result, and this proves the remaining assertions
of the theorem.
 
{\bf Proof of Theorem \ref{theorem 3.21.01}}.
After \eqref{6.14.1} has been proved and the modification
$u_{t}(x-\sigma b_{t})$ of $\pi_{t}(x)$ has been chosen
the assertions related to \eqref{3.24.1} and \eqref{3.24.02} 
follow directly from \eqref{6.12.1} and \eqref{6.12.2}.
This proves the theorem.

\mysection{Appendix}
                                    \label{section 6.14.3}

Let $w^{1}_{t},w^{2}_{t},...$ be independent one-dimensional
Wiener processes with respect to $\cF_{t}$, and let
$\cP$ denote the predictable $\sigma$-field related to
$\{\cF_{t}\}$. Assume that  on $\Omega\times
(0,\infty)\times\bR^{d}$ we are given
$a_{t}(x)=(a^{ij}_{t}(x))$ which is a
 $d\times d$-symmetric matrix  valued
function, $b_{t}(x)$ which is an $\bR^{d}$-valued
function, real-valued $c_{t}(x)$, and 
$\sigma^{i\cdot}_{t}(x)
=(\sigma^{i,1}_{t}(x),\sigma^{i,2}_{t}(x),...)$, $i=1,...,d$,
and $\nu_{t}(x)=(\nu^{1}_{t}(x),\nu^{2}_{t}(x),...)$,
which are $\ell_{2}$-valued functions.

\begin{assumption}
                                      \label{assumption 5.27.1}
(i) The functions $a,b,c,\sigma,\nu$
 are bounded and measurable
as functions of $(\omega,t,x)$ and are 
predictable  
as functions of $ (\omega,t) $ for each $x$.

(ii) The functions $a_{t}(x)$ are Lipschitz
continuous in $x$ with constant $K$ (independent of
$\omega,t$).

(ii) For any $\lambda,x\in\bR^{d}$, $t\geq0$,
and $\omega$
\begin{equation}
                                      \label{5.26.1}
(2a^{ij}_{t}-\alpha^{ij}_{t})(x)\lambda^{i}\lambda^{j}
\geq \delta_{0}|\lambda|^{2},
\end{equation}
where $\delta_{0}$ is a constant, $\delta_{0}>0$,
$\alpha^{ij}:=\sigma^{ik}\sigma^{jk}$.

\end{assumption}

Let $G\subset\bR^{d}$ be an open set, $T\in[0,\infty)$, and
let $u_{t}(x)$ be a real-valued function
given on $\Omega\times[0,T]\times\bar{G}$.
Also suppose that on $\Omega\times(0,T)\times G$
we are given
 functions $f_{t}(x)$ and
$g_{t}(x)=(g^{1}_{t}(x),g^{2}_{t}(x),...)$
with values
in $\bR$ and $ \ell_{2}$, respectively.

Define $D_{i}=\partial/\partial x^{i}$, $D_{ij}=D_{i}D_{j}$,
$Du$ the gradient of $u$, $D^{2}u$ its Hessian.
   
\begin{assumption}
                                      \label{assumption 5.27.2}
(i) For any $\omega$ the function $u$
is continuous in $[0,T]\times\bar{G}$ and  vanishes
on $[0,T]\times\partial G$ and, moreover (if $G$ is unbounded),
for any $\delta>0$ and $\omega$ there exists a compact set $\Gamma
\subset G$ such that $|u_{t}(x)|\leq\delta$
for all $t\in[0,T]$ and $x\in G\setminus\Gamma$.
For any $x$
the function $u_{t}(x)$ is $\cF_{t}$-adapted.
For any $\omega$ (if $G$ is unbounded)
$$
\int_{G}|u_{0}|^{2}\,dx+\int_{0}^{T}\int_{G}|u_{t}|^{2}\,dxdt<\infty.
$$

(ii) For any $\omega$, for almost any $t\in[0,T]$,
 the second-order derivatives
$D^{2}u_{t}(x)$ are continuous with respect to $x$ and
for any compact set $\Gamma\subset G$
$$
\int_{0}^{T}\int_{\Gamma}|D^{2}u_{t}|\,dxdt<\infty.
$$

(iii) For   each   $\omega$ and compact set $\Gamma\subset G$
$$
\int_{0}^{T}\int_{\Gamma}|Du_{t}|^{2}\,dx dt<\infty.
$$

(iii) The functions $f_{t}(x)$ and $g_{t}(x)$ are 
$\cP\otimes\cB(G)$-measurable
as functions of $(\omega,t,x)$, where 
$\cB(G)$ is the Borel $\sigma$-field on $G$.
For each  
 $x\in G$  and $\omega$ we have
$$ 
\int_{0}^{T}(|D^{2}u_{s}(x)|+|Du_{s}(x)|^{2}+|f_{s}(x)|^{2} +
|g_{s}(x)|^{2}_{ \ell_{2}})\,ds<\infty.
$$

(iv) For each  $\omega$
$$
\int_{G}\int_{0}^{T}(|f_{s}(x)|^{2} +
|g_{s}(x)|^{2}_{ \ell_{2}})\,dsdx<\infty,
$$
and for each $x\in G$ with probability one  we have
for all $t\in [0,T]$ that
$$
u_{t}(x)=u_{0}(x)
+\int_{0}^{t}\big[\sigma^{ik}_{s}(x)D_{i}u_{s}(x)
+\nu^{k}_{s}(x)u_{s}(x)+g^{k}_{s}(x)\big]\,dw^{k}_{s}
$$
$$
+\int_{0}^{t}\big[a^{ij}_{s}(x)
D_{ij}u_{s}(x)+b^{i}_{s}(x)D_{i}u_{s}(x)
+c_{s}(x)u_{s}(x)+f_{s}(x)\big]\,ds.
$$

\end{assumption}

\begin{theorem}
                                       \label{theorem 5.28.1}
Under the above assumptions
$$
\int_{0}^{T}\int_{G}|Du_{t}(x)|^{2}\,dxdt<\infty
$$
(a.s.) and for any $\phi\in C^{\infty}_{0}(G)$
with probability one
$$
(u_{t} ,\phi)=(u_{0} ,\phi)
+\int_{0}^{t}(\phi,
\sigma^{ik}_{s}D_{i}u_{s} 
+\nu^{k}_{s}u_{s} +g^{k }_{s}) \,dw^{k}_{s}
$$
$$
+\int_{0}^{t}\big[(\phi,
(b^{i}_{s}-D_{j}a^{ij} ) D_{i}u _{s} 
+c_{s}u _{s}
 +f _{s})-(D_{j}\phi,a^{ij}_{s}D_{i}
u_{s} )\big]\,ds,
$$
for all $t\in[0,T]$.
\end{theorem}

Proof.
By the stochastic Fubini theorem (see, 
for instance, Lemma 2.7  of \cite{Kr11}),
for any $\phi\in C^{\infty}_{0}(G)$ 
with
probability one for all $t\geq 0$
$$
(u_{t}^{\delta},\phi)=(u_{0}^{\delta},\phi)
+\int_{0}^{t}(\phi,
\sigma^{ik}_{s}D_{i}u_{s} 
+\nu^{k}_{s}u_{s} +g^{k }_{s}) \,dw^{k}_{s}
$$
$$
+\int_{0}^{t}\big[(\phi,
(b^{i}_{s}-D_{j}a^{ij} ) D_{i}u _{s} 
+c_{s}u _{s}
 +f _{s})-(D_{j}\phi,a^{ij}_{s}D_{i}
u_{s} )\big]\,ds,
$$
where $u^{\delta}=u-\delta$.

For $\varepsilon,\delta>0$
define 
$$
K_{\varepsilon}=\{x\in G:\dist(x,\partial D)\leq
\varepsilon\quad\text{or}\quad|x|\geq\varepsilon^{-1}\},
$$
$$
\tau_{\varepsilon,\delta}
=T\wedge\inf\{t\geq0:\sup_{G\setminus K_{\varepsilon}}
|u_{t}|\geq\delta\}.
$$
Observe that as $\varepsilon\downarrow0$ we have 
$$
\sup_{G\setminus K_{\varepsilon}}
|u_{t}|\downarrow0
$$
uniformly with respect to $t\in[0,T]$ by assumption.
Therefore $\tau_{\varepsilon,\delta}\to T$ 
as $\varepsilon\downarrow0$ for any
$\delta>0$. Also notice that
$$
|u_{t}(x)|\leq \delta
$$
in $G\setminus K_{\varepsilon}$ if 
$0\leq t<\tau_{\varepsilon,\delta} $.

By Lemma 2.5  of \cite{Kr07} for any 
$\phi\in C^{\infty}_{0}(G)$
(a.s.) for all $t\in[0,T]$
\begin{equation}
                                                 \label{5.28.3}
\|(\phi u^{\delta}
_{t\wedge\tau_{\varepsilon,\delta}})^{+}\|_{L_{2}}^{2}=
\|(\phi u_{0}^{\delta})^{+}\|_{L_{2}}^{2}+
\int_{0}^{t}h _{s}\,ds+m _{t} ,
\end{equation}
where
$$
h _{s}:=2\big((\phi u^{\delta}_{s})^{+}, \phi\{
(b^{i}_{s}-D_{j}a^{ij} ) D_{i}u _{s} 
+c_{t}u _{s}
 +f _{s}\}
$$
$$
-a^{ij}_{t}(D_{i}
u_{s} )D_{j}\phi\big)I_{s<\tau_{\varepsilon,\delta}}
-2(I_{\phi u^{\delta}_{s}>0}
D_{i}(\phi u^{\delta}_{s}),\phi 
a^{ij}_{s}D_{j}
u_{s} )
I_{s<\tau_{\varepsilon,\delta}}
$$
$$
+\sum_{k}\big\|
\phi\big[
\sigma^{ik}_{s}D_{i}u _{s} 
+\nu^{k}_{s}u_{s} +g^{k }_{s})
\big]I_{\phi u^{\delta}_{s}>0}\big\|_{L_{2}}^{2}
I_{s<\tau_{\varepsilon,\delta}} ,
$$
$$
m_{t} :=2\int_{0}^{t}I_{s<\tau_{\varepsilon,\delta}}
(\phi(\phi u^{\delta})^{+}_{s}, 
\sigma^{ik}_{s}D_{i}u _{s} 
+\nu^{k}_{s}u_{s} 
+g^{k }_{s})\,dw^{k}_{s}.
$$
We take $\phi\geq0$ such that $\phi=1$ on 
$K_{\varepsilon}$. Then for 
$0\leq s<\tau_{\varepsilon,\delta}$
$$
(\phi u^{\delta}_{s})^{+}= u^{\delta}_{s} I_{u_{s}>\delta}.
$$
Indeed, if $u^{\delta}_{s}(x)\leq0$, then both sides
vanish. However, if $u^{\delta}_{s}(x)>0$ and 
$0\leq s<\tau_{\varepsilon,\delta}$, then 
$x\in K_{\varepsilon}$   and $\phi(x)=1$.
By also taking into account that
$$
I_{u _{s}>\delta,s<\tau_{\varepsilon,\delta}}
D_{j}\phi=0 
$$ 
we transform \eqref{5.28.3} for such a $\phi$ into
\begin{equation}
                                                 \label{5.28.4}
\|(u^{\delta}
_{t\wedge\tau_{\varepsilon,\delta}})^{+}\|_{L_{2}}^{2}=
\|( u_{0}^{\delta})^{+}\|_{L_{2}}^{2}+
\int_{0}^{t}h _{s}\,ds+m _{t} ,
\end{equation} 
where
$$
h _{s}:=2\big( u^{\delta}_{s}I_{u_{s}>\delta},  
(b^{i}_{s}-D_{j}a^{ij}_{ij}) I_{u_{s}>\delta} D_{i}
 u _{s}   
+c_{t}u _{s}
 +f _{s}\big)I_{s<\tau_{\varepsilon,\delta}}
$$
$$
-2( I_{u_{s}>\delta} D_{i}
 u _{s} ,  
a^{ij}_{t} D_{j}
 u _{s}  )
I_{s<\tau_{\varepsilon,\delta}}
$$
$$
+\sum_{k}\big\|
\big[
\sigma^{ik}_{s}D_{i} u^{\delta}_{s}  
+\nu^{k}_{s}u_{s} +g^{k }_{s})
\big]I_{u _{s}>\delta}\big\|_{L_{2}}^{2}
I_{s<\tau_{\varepsilon,\delta}} ,
$$
$$
m_{t} :=2\int_{0}^{t}I_{s<\tau_{\varepsilon,\delta}}
(( u^{\delta})^{+}_{s}, 
\sigma^{ik}_{s}D_{i}u _{s} 
+\nu^{k}_{s}u_{s}^{\delta}
+g^{k\delta}_{s})\,dw^{k}_{s}.
$$

To proceed further we recall again the fact that
the two conditions: 
$s<\tau_{\varepsilon,\delta}$ and $u_{s}(x) >\delta$,
imply that $x\in K_{\varepsilon}$ so that
$$
\int_{0}^{\tau_{\varepsilon,\delta}}\int_{G}
|D u _{t} |^{2}I_{u_{t}>\delta}\,dxdt
\leq
\int_{0}^{T}\int_{K_{\varepsilon}}|Du_{t}|^{2}\,dxdt<\infty.
$$
After that we use \eqref{5.26.1} and the inequalities like
$ab\leq\varepsilon a^{2}+\varepsilon^{-1}b^{2}$.
We also recall that $b^{i}_{s}-D_{j}a^{ij}_{s}, c_{s}$,
and $|\nu_{s}|_{\ell_{2}}$ are bounded and then in an absolutely
standard way derive from \eqref{5.28.4} that
there exists a constant $N\in(0,\infty)$, 
independent of $\omega,t,
\varepsilon,\delta$, such that for any $t\in[0,T]$ and $\omega$
$$
N^{-1}\int_{0}^{t\wedge\tau_{\varepsilon,\delta}}
\|I_{u_{s}>\delta}Du_{s}\|_{L_{2}}^{2}\,ds
\leq N\| u_{0} \|_{L_{2}}^{2}+
N\int_{0}^{t\wedge\tau_{\varepsilon,\delta}}
\| u_{s}\|_{L_{2}}^{2}\,ds+m_{t}.
$$ 
In particular, for any stopping time $\tau\leq T$
$$
N^{-1}\int_{0}^{\tau\wedge\tau_{\varepsilon,\delta}}
\|I_{u_{s}>\delta}Du_{s}\|_{L_{2}}^{2}\,ds
\leq N\| u_{0} \|_{L_{2}}^{2}I_{\tau>0}+
N\int_{0}^{\tau\wedge\tau_{\varepsilon,\delta}}
\| u_{s}\|_{L_{2}}^{2}\,ds+m_{\tau}.
$$  
If $\tau$ is a localizing time for the
local martingale $m_{t}$ starting at zero,
then $Em_{\tau}=0$ and
$$
N^{-1}E\int_{0}^{\tau\wedge\tau_{\varepsilon,\delta}}
\|I_{u_{s}>\delta}Du_{s}\|_{L_{2}}^{2}\,ds
\leq NE\| u_{0} \|_{L_{2}}^{2}I_{\tau>0}+
NE\int_{0}^{\tau\wedge\tau_{\varepsilon,\delta}}
\| u_{s}\|_{L_{2}}^{2}\,ds.
$$
The last inequality, actually,
 holds for any stopping time $\tau\leq T$,
which is easily proved by approximation.

Now we first let $\varepsilon\downarrow0$ and then
$\delta\downarrow0$. Then by the monotone convergence theorem
we obtain
$$
N^{-1}E\int_{0}^{\tau }
\|I_{u_{s}>0}Du_{s}\|_{L_{2}}^{2}\,ds
\leq NE\| u_{0} \|_{L_{2}}^{2}I_{\tau>0}+
NE\int_{0}^{\tau }
\| u_{s}\|_{L_{2}}^{2}\,ds.
$$
Similarly,
$$
N^{-1}E\int_{0}^{\tau }
\|I_{u_{s}<0}Du_{s}\|_{L_{2}}^{2}\,ds
\leq NE\| u_{0} \|_{L_{2}}^{2}I_{\tau>0}+
NE\int_{0}^{\tau }
\| u_{s}\|_{L_{2}}^{2}\,ds,
$$
and since one  knows that $I_{u=0}Du=0$ (a.e.)
for any function $u$ of class $W^{1}_{2}$ we finally conclude that
\begin{equation}
                                              \label{5.29.1}
N^{-1}E\int_{0}^{\tau }
\| Du_{s}\|_{L_{2}}^{2}\,ds
\leq NE\| u_{0} \|_{L_{2}}^{2}I_{\tau>0}+
NE\int_{0}^{\tau }
\| u_{s}\|_{L_{2}}^{2}\,ds.
\end{equation}
For $n>0$ and $\tau=\tau^{n}$, where
$$
\tau_{n}=T\wedge\inf\{t\geq0:\|u_{0}\|_{L_{2}}+
\int_{0}^{t}\|u_{s}\|_{L_{2}}^{2}\,ds\geq n\},
$$
the right-hand side of \eqref{5.29.1} is finite. Hence
$$
\int_{0}^{\tau^{n} }
\| Du_{s}\|_{L_{2}}^{2}\,ds<\infty
$$
(a.s.). To prove the first assertion
of the theorem, it only remains to observe that,
 by assumption, for any $\omega$, we have
 $\tau_{n}=T$ for all sufficiently large $n$.
The second assertion follows from the stochastic Fubini theorem.
The theorem is proved.

\end{document}